\documentclass [12pt,leqno] {article}

\usepackage{amsfonts,url}
\def\supp{\mathop{\rm supp}\nolimits}

\newcommand{\qed}{~\hfill~$\fbox{}$ }

\newcommand{\proces}{( X_t,P^x)} 

\newcommand{\R}{ \mathbb{R}^{d}}
\newcommand{\N}{\mathbb{N}}

\newcommand{\dowod}{{\em Proof}.\/ }
\newcommand{\indyk}[1]{{\bf 1}_{#1}}
\newcommand{\sfera}{ \mathbb{S}}

\newcommand{\Fourier}{ {\cal F}}

\newcommand{\nubounded}[1]{\bar{\nu}_{#1}}

\newcommand{\scalp}[2]{#1\cdot#2}
\newcommand{\qd}{\underline{q}}
\newcommand{\qu}{\bar{q}}
\newcommand{\czc}[1]{\left\lfloor {#1} \right\rfloor}

\newtheorem{lemat}{\indent\sc Lemma}

\newtheorem{twierdzenie}{\indent\sc Theorem}
\newtheorem{wniosek}[lemat]{\indent\sc Corollary}

\newcounter{conum} 

\begin{document}

\title{Estimates of tempered stable densities}
\author{Pawe{\l} Sztonyk}
\footnotetext{ Institute for Mathematical Stochastics, TU Dresden,
  and\\ 
  Institute of Mathematics and Computer Science,
  Wroc{\l}aw University of Technology,
  Wybrze{\.z}e Wyspia{\'n}\-skie\-go 27,
  50-370 Wroc{\l}aw, Poland\\
  {\rm e-mail: sztonyk@pwr.wroc.pl}
}
% \footnotetext{ Institute of Mathematics,
%   Wroc{\l}aw University of Technology,
%   Wybrze{\.z}e Wyspia{\'n}\-skie\-go 27,
%   50-370 Wroc{\l}aw, Poland\\
%   {\rm e-mail: sztonyk@im.pwr.wroc.pl,
%     bogdan@im.pwr.wroc.pl}}
\date{March 18, 2008}
\maketitle

\begin{center}
  Abstract
\end{center}
\begin{scriptsize}
  Estimates of densities of convolution semigroups of probability measures are given under specific assumptions
  on the corresponding L\'evy measure and the L\'evy--Khinchin exponent. The assumptions are satisfied, e.g., by
  tempered stable semigroups of J. Rosi{\'n}ski.
\end{scriptsize}

\footnotetext{2000 {\it MS Classification}:
Primary 60G51, 60E07; Secondary 60J35, 47D03, 60J45 .\\
%47D03 Groups and semigroups of linear operators
%31C05 Harmonic, subharmonic, superharmonic functions
%60J35 Transition functions, generators and resolvents 
%60G51 Processes with independent increments
{\it Key words and phrases}: stable process, tempered stable process, semigroup of measures, transition density.\\
The research was supported by The State Committee for Scientific Research (Poland, KBN 1 P03A 026 29) and The Alexander von Humboldt Foundation (Germany, 3-PL/1122470).
}
%----------------------------------------------------------------

\section{Introduction and main results}\label{Prel}

Let $\nu$ be a symmetric L\'evy measure on $\R$ where $d\in\{1,2,3,\dots\}$.
%We will assume that $\nu$ is positive and symmetric:
%$$
%\nu(D)=\nu(-D)\,,\quad D\subset {\cal B}(\R)\,.
%$$

We consider the convolution semigroup of probability measures $\{ P_t, t\geq 0 \}$ with the Fourier transform $\Fourier(P_t)(\xi)=\exp(-t\Phi(\xi))$, where
%&  =  & -\int \left(e^{i\scalp{\xi}{y}}-1-i\scalp{\xi}{y}\indyk{B(0,1)}(y)\right)\nu(dy) \\
$$
  \Phi(\xi) =   \int \left(1-\cos(\scalp{\xi}{y})\right)\nu(dy).
$$
%Here $\scalp{\xi}{y}$ denotes the usual inner product of $\xi,y\in\R$. 
The semigroup determines 
the stochastic L\'evy process $\proces$ on $\R$ with the generating triplet $(0,\nu,0)$ (we use
the terminology of \cite{Sato})
and transition probabilities $P(t,x,A)=P_t(A-x)$. If
$$
  \nu(D)= \int_{\sfera} \int_0^\infty \indyk{D}(s\theta) s^{-1-\alpha}ds\mu(d\theta),
$$
where $\alpha\in (0,2)$ and $\mu$ is a bounded measure on $\sfera=\{x\in\R\: : |x|=1\}$,
then $\proces$ is called the {\it $\alpha$--stable L\'evy process.} If $\mu$ is nondegenerate, i.e., if there is no proper
linear subspace $M$ of $\R$ such that $\supp(\mu)\subset M$, then the $\alpha$--stable
symmetric measure $P_t$ is absolutely continuous for every $t>0$ and the corresponding density
$p_t$ is smooth, bounded and it has the {\it scaling property}:
$p_t(x)=t^{-d/\alpha}p_1(t^{-1/\alpha}x)$.
 
Stable processes are important tool in theoretical probability, physics and finance and their asymptotic
properties are subject of interest of many papers.
W.E. Pruitt and S.J. Taylor investigated in \cite{PrTa} multivariate stable densities in the general setting. They obtained the estimate $p_1(x)\leq c (1+|x|)^{-1-\alpha}$ by Fourier-analytic methods. Indeed such a decay can be obtained if the spectral measure $\mu$ has an atom (see the estimates from below in \cite{Hi1} and \cite{Hi2}). In the well--known case of the rotation invariant $\alpha$--stable process with uniform $\mu$ we however have
$p_1(x)\approx (1+|x|)^{-d-\alpha}$ (see \cite{Zol}). Using the perturbation formula P. G{\l}owacki and W. Hebisch proved in \cite{G} and \cite{GH} that if $\mu$
has a bounded density, $g_{\mu}$, with respect to the surface measure on $\sfera$ then $p_1(x)\leq c (1+|x|)^{-d-\alpha}$. When $g_{\mu}$ is continuous on $\sfera$ we even have
$\lim_{r\to\infty}r^{d+\alpha}p_1(r\theta)=cg_{\mu}(\theta)$, $\theta\in\sfera$ and if $g_{\mu}(\theta)=0$ then additionally $\lim_{r\to\infty}r^{d+2\alpha}p_1(r\theta)=c_{\theta}>0$, which was proved by J. Dziuba\'nski in \cite{Dziub}.
A. Zaigraev in \cite{Zai} obtained further asymptotic expansions of the $\alpha$--stable density for sufficiently regular $g_\mu$.

More recent asymptotic results for stable L\'evy processes are given 
in papers \cite{W} and \cite{BS2007}. In particular if for some $\gamma\in[1,d]$ 
the measure $\nu$ is a $\gamma$--measure on $\sfera$ , i.e., 
$$
  \nu(B(x,r))\leq c r^{\gamma} \quad  \mbox{for every}\quad x \in\sfera,\, r\leq 1/2,
$$  
or equivalently 
$$
\mu(B(\theta,r)\cap \sfera)\leq c r^{\gamma-1},\quad  \theta\in\sfera,\, r\leq 1/2,
$$ 
then we have 
$$
  p_1(x)\leq c (1+|x|)^{-\alpha-\gamma},\quad x\in\R.
$$ 
By scaling $p_t(x)\leq c t^{-d/\alpha}(1+t^{-1/\alpha}|x|)^{-\alpha-\gamma}$ for every $t>0$. It follows also from 
\cite[Theorem 1.1]{W} that if for some $\theta_0\in\sfera$ we have
$$
\mu(B(\theta_0,r)\cap \sfera)\geq c r^{\gamma-1},\quad  r\leq 1/2,
$$
then 
$$
  p_1(r\theta_0)\geq c (1+r)^{-\alpha-\gamma},\quad r>0.
$$

Our main goal is to extend the estimates to more general class
of semigroups and processes. The obtained below results
cover a wide class of examples which we discuss in more detail in the Section \ref{sec:Temp}.

We fix the constant $\alpha\in (0,2)$ and we always assume that 
there exists a positive $c$ such that
\begin{equation}\label{eq:Ass0}
  \Phi(\xi)\geq c |\xi|^{\alpha}\,\, \mbox{for}\,\, |\xi|>1.
\end{equation}
Note that
(\ref{eq:Ass0}) is satisfied, e.g., if we have 
$$
 \int_{|y|\leq r}|\scalp{\theta}{y}|^2\nu(dy)\geq c r^{2-\alpha}\quad  \mbox{for}\quad 
  r \leq 1,\, \theta\in\sfera.
$$
It follows from (\ref{eq:Ass0}) that the measures $P_t$ are absolutely continuous with respect to the Lebesgue measure and their densities $p_t$ are smooth and bounded.

All the sets, functions and measures considered in the sequel will be Borel.

\begin{twierdzenie}\label{th:PEstimates}
Let there exist constants $c>0$, $\gamma\in [1,d]$, $K\in[1,\infty)$, a~bounded measure $\mu$
on $\sfera$, and
a~bounded nonincreasing function $\qu:(0,\infty)\to (0,\infty)$, such that
\begin{equation}\label{eq:Ass1}
  \nu(A) \leq \int_{\sfera}\int_0^\infty \indyk{A}(s\theta) s^{-1-\alpha}\qu(s)ds\mu(d\theta),
  \quad A\subset\R,
\end{equation}
\begin{equation}\label{eq:Ass2}
  \qu(s) \leq K\qu(2s),\quad s>0,
\end{equation}
and
\begin{equation}\label{eq:Ass3}
  \mu(B(\theta,r)\cap\sfera) \leq c r^{\gamma-1},\quad \theta\in\sfera,\, r<1/2.
\end{equation}
Then there exists a constant $C$ such that
\begin{eqnarray}\label{eq:AbDensityEstimate}
  p_t(x) & \leq    & C \min\left(t^{-d/\alpha},t^{1+\frac{\gamma-d}{\alpha}}
                     |x|^{-\alpha-\gamma}\,\qu(|x|)\right)\nonumber \\
         & \approx & t^{-d/\alpha}\left(1+t^{-1/\alpha}|x|\right)^{-\alpha-\gamma}\qu(|x|), 
                     \quad x\in\R,\, t\in(0,1).
\end{eqnarray}
\end{twierdzenie}

The {\it doubling property} (\ref{eq:Ass2}) is equivalent to the following:
there exist $c>0$ and $\eta\geq 0$ such that
\begin{equation}\label{eq:Ass2a}
   \frac{\qu(r)}{\qu(R)} \leq c\left(\frac{r}{R}\right)^{-\eta},\quad 0<r\leq R.
\end{equation}
The typical examples of $\qu$ are $\qu(s)=(1+s)^{-a},$ or $\qu(s)=(\log(e+s))^{a}(1+s)^{-ma},$ for $a\geq 0$ and
$m> 1$.

We note that if (\ref{eq:Ass1}), (\ref{eq:Ass2}) and (\ref{eq:Ass3}) hold then we have
\begin{equation}\label{eq:Col1}
  \nu(B(x,r)) \leq c r^\gamma |x|^{-\alpha-\gamma}\qu(|x|),\quad r<\frac{1}{2}|x|,\,
  x\in\R\setminus\{0\},
\end{equation}
and
\begin{equation}\label{eq:Col2}
  \nu(B(0,r)^c)\leq c r^{-\alpha}\qu(r), \quad r\in(0,\infty).
\end{equation}
We omit the easy proof. The partial converse of Theorem \ref{th:PEstimates} is given in the following theorem.

\begin{twierdzenie}\label{th:PEstimatesB}
Let there exist a set $A\subset\R$, constants $\gamma\in[1,d]$, $c_1,c_2>0$ and 
a bounded function $\qd :(0,\infty)\to (0,\infty)$ such that
\begin{equation}\label{eq:Ass6}
  \nu(B(x,r)) \geq c_1 r^\gamma |x|^{-\alpha-\gamma}\qd(|x|),\quad x\in A,\, r>0,
\end{equation}
\begin{equation}\label{eq:Ass7}
  \nu(B(0,r)^c)\leq c_2 r^{-\alpha}, \quad r\in (0,1).
\end{equation}
%and
%If the assumptions (\ref{eq:Ass6}), (\ref{eq:Ass4a}) and (\ref{eq:Ass5old}) are satisfied 
Then there exists a constant $C$ such that
\begin{equation}\label{eq:BeDensityEstimate}
  p_t(x)    \geq   C \min\left(t^{-d/\alpha},t^{1+\frac{\gamma-d}{\alpha}}
                     |x|^{-\alpha-\gamma}\,\qd(|x|)\right), 
                     \quad  x\in A,\, t\in(0,1).
\end{equation}
\end{twierdzenie}

The mild assumptions on the function $\qd$ allow to use Theorem \ref{th:PEstimatesB}
for a wide class of processes. 
An important example is 
the relativistic $\alpha$--stable process (see \cite{Ryz,GrzRyz08}) with $\Phi(\xi)=(|\xi|^2+1)^{\alpha/2}-1$. It follows from \cite{GrzRyz08} that
in this case (\ref{eq:Ass6}) and (\ref{eq:Ass7}) hold with $\gamma=d$, $A=\R$, and $\qd(s)=(1+s)^{\frac{d+\alpha-1}{2}}e^{-2s}$, and
we get $ p_t(x) \geq c \min\left(t^{-d/\alpha},t|x|^{-\alpha-d}(1+|x|)^{\frac{d+\alpha-1}{2}}e^{-2|x|} \right)$
for $x\in\R$ and $t\in (0,1)$. We discuss this example in detail in Section \ref{sec:Temp}.

The above theorems hold for small times $t$ and 
%in our general case we can not use the scaling property
%to extend the results to all values of $t$. Therefore we
%obtain the estimates in two steps considering 
below we consider the case of large $t$. We assume still (\ref{eq:Ass0}) which guarantees in particular the existence and smoothness of densities. However the behaviour of $p_t$ at the origin 
for large $t$ depends on 
the asymptotic of $\Phi(\xi)$ for small $\xi$ (see Lemma \ref{lm:pBasicEstimateBigt}) where $\Phi$ may decay faster
then $|\xi|^\alpha$. Hence we strengthen below our assumptions.

\begin{twierdzenie}\label{th:PEstimatesBigt}
Let (\ref{eq:Ass1}), (\ref{eq:Ass2}) and (\ref{eq:Ass3}) hold and let there exist constants $c$ and $\beta\in[\alpha,2]$ such that
\begin{equation}\label{eq:Ass5} 
  \int_1^\infty s^{\beta-\alpha-1}\qu(s)<\infty,
\end{equation}
and
\begin{equation}\label{eq:Ass4}
  \Phi(\xi)\geq c |\xi|^{\beta},\quad |\xi|\leq 1.
\end{equation}
%If the assumptions (\ref{eq:Ass1old}),(\ref{eq:Ass2old}), (\ref{eq:Ass0})  
%and (\ref{eq:Ass2}) are satisfied 
Then there exists a constant $C$ such that
\begin{equation}\label{eq:AbDensityEstimateBigt}
  p_t(x)\leq C \min\left(t^{-d/\beta},t^{1+\frac{\gamma-d}{\beta}}|x|^{-\alpha-\gamma}\qu(|x|)\right), \quad
  x\in\R,\, t\in(1,\infty).
\end{equation}
\end{twierdzenie}

\begin{twierdzenie}\label{th:PEstimatesBBigt}
Let (\ref{eq:Ass6}) hold and let there exist constant $\beta\in[\alpha,2]$ and c such that
\begin{equation}\label{eq:Ass9}
  \nu(B(0,r)^c)\leq c r^{-\beta}, \quad r\in (0,\infty),
\end{equation}
\begin{equation}\label{eq:Ass10}
  c^{-1}|\xi|^\beta \leq \Phi(\xi) \leq c |\xi|^{\beta}, \quad |\xi|\leq 1.
\end{equation}
%If the assumptions (\ref{eq:Ass6}), (\ref{eq:Ass4a}) and (\ref{eq:Ass5old}) are satisfied 
Then there exists a constant $C$ such that
\begin{equation}\label{eq:BeDensityEstimateBigt}
  p_t(x)\geq C \min(t^{-d/\beta},t^{1+\frac{\gamma-d}{\beta}}|x|^{-\alpha-\gamma}\qd(|x|)), \quad
  x\in A,\, t\in(1,\infty).
\end{equation}
\end{twierdzenie}

We like to mention related recent results. Z.-Q.~Chen and T.~Kumagai in \cite{ChKum} and \cite{ChKum08} investigate the case of symmetric jump--type Markov processes 
on some class of metric measure spaces with jump kernels which are not translation invariant. They obtain estimates
of the densities in small time (see Theorem~1.2 in \cite{ChKum08}) which are similar to given above. However the jump kernels in \cite{ChKum} and \cite{ChKum08}  
are assumed to have densities (which corresponds to $\gamma=d$ in our setting) comparable with certain rotation invariant functions. Here the assumptions (\ref{eq:Ass1}) and (\ref{eq:Ass3}) include the isotropic estimates only from above and in Theorems
\ref{th:PEstimatesB} and \ref{th:PEstimatesBBigt} the estimate (\ref{eq:Ass6}) can depend on a chosen direction so we allow for more anisotropy here. The results for large $t$ are restricted in $\cite{ChKum08}$ to the case of jump densities which decay in infinity not faster then $|x|^{-2-d}$.

The asymptotics of the densities of the truncated stable--like processes have been investigated in \cite{ChKK} by Z.-Q.~Chen, P.~Kim and T.~Kumagai. K.~Bogdan and
T.~Jakubowski in \cite{BJ2007} obtained estimates of the heat kernels of the fractional Laplacian perturbed by gradient operators. The derivatives of the stable densities have been considered in \cite{Lewand} and \cite{Szt2007}.

J.~Picard \cite{Pi} studied the density in small time for jump processes obtaining the estimates
for solutions of certain stochastic differential equation driven by stable L\'evy process. Here we improve the methods used in \cite{Pi} and \cite{BS2007}
and extend the results to a wide class of L\'evy non--stable processes. 
Our assumptions are satisfied, e.g., by the tempered stable processes of J. Rosi\'nski (\cite{Ros07}) and the layered
stable processes of C.~Houdr\'e and R.~Kawai (\cite{HoKawa}).

The paper is organized as follows. In Section \ref{sc:Local} we investigate densities
of the semigroup determinated by the truncated L\'evy measure $\indyk{B(0,\varepsilon)}\nu$ where
$\varepsilon$ is comparable with $t^{1/\alpha}$ for small $t$ and $t^{1/\beta}$ for large $t$. In Section \ref{sc:BP}
we consider $\indyk{B(0,\varepsilon)^c}\nu$ which is a bounded measure and therefore the corresponding semigroup 
is given by an exponential formula.
We collect these results in Section \ref{sc:MainR} and prove the main theorems. Examples
are given in Section \ref{sec:Temp}.

We use $c,C$ (with subscripts) to denote finite positive constants
which depend only on the measure $\mu$, the functions $\qu,\qd$, constants $\gamma,\alpha,\beta$, and the dimension $d$. Any {\it additional} dependence
is explicitly indicated by writing, e.g., $c=c(n)$.
The value of $c$, when used without subscripts, may
change from place to place. We write $f\approx g$ to indicate that there is a $c$
such that $c^{-1}f \leq g \leq c f$.

%--------------------------------------------------------------
\section{Local part}\label{sc:Local}
For $\varepsilon>0$ we denote $\tilde{\nu}_\varepsilon=\indyk{B(0,\varepsilon)}\nu$. 
We consider the corresponding semigroup of measures $\{\tilde{P}^\varepsilon_t,\; t\geq 0\}$ such that
\begin{equation}\label{eq:FPtilde}
  \Fourier(\tilde{P}^\varepsilon_t)(\xi) =
  \exp\left(t \int (\cos(\scalp{\xi}{y})-1)
  \tilde{\nu}_\varepsilon(dy)\right)\, ,
  \quad \xi\in\R\, .
\end{equation}

By (\ref{eq:Ass0}) we get
\begin{eqnarray}\label{eq:FTildePEstimate}
  \Fourier(\tilde{P}^\varepsilon_t)(\xi) 
  &   =  & \exp\left(-t\int_{|y|<\varepsilon}
           (1-\cos(\scalp{y}{\xi}))
           \nu(dy)\right) \nonumber \\
  &   =  & \exp\left(-t\left(\Phi(\xi)-\int_{|y|\geq\varepsilon}
           (1-\cos(\scalp{y}{\xi}))
           \nu(dy)\right)\right) \nonumber \\
  & \leq & \exp(-t\Phi(\xi))\exp(2t\nu(B(0,\varepsilon)^c)) \nonumber \\
  & \leq & \exp(-ct|\xi|^\alpha)\exp(2t\nu(B(0,\varepsilon)^c)),%e^{ct\varepsilon^{-\alpha}}
           \quad |\xi|>1.
\end{eqnarray}
From (\ref{eq:FTildePEstimate}), 
\cite[Theorem 3.7.13]{Jc1} and the multivariable version of Faa di Bruno's formula (see, e.g., \cite{ConSav96}) it follows  that for every $\varepsilon>0$ and $t>0$ the measure $\tilde{P}^\varepsilon_t$ is absolutely continuous with respect to
the Lebesgue measure and its density $\tilde{p}^\varepsilon_t$ belongs to the Schwarz space of smooth
rapidly decreasing functions.

We will often use $\tilde{P}^\varepsilon_t$ and $\tilde{p}^\varepsilon_t$ with $\varepsilon=t^{1/\alpha}$.
For simplification we will denote 
$$ 
  \tilde{P}_t=\tilde{P}^{t^{1/\alpha}}_t\,\, \mbox{and}\,\, \tilde{p}_t=\tilde{p}^{t^{1/\alpha}}_t.
$$

Since $\tilde{p}_t$ belongs to the Schwarz space we certainly have that $\tilde{p}_t(y)\leq c(t,m)(1+|y|)^{-m}$ for every
$m>0$ and $t>0$.
We improve the estimates in the following two Lemmas.

\begin{lemat}\label{lm:MomentsTildeP}
If we have
\begin{equation}\label{eq:AssAux1}
  \nu(B(0,r)^c) \leq c r^{-\alpha},\,\, \mbox{for}\,\, r\in(0,\infty),
\end{equation}
then for every $n\in\N$ there exists a constant $c_n=c(n)$ such that
\begin{equation}\label{eq:MomentsTildeP}
  \int_{\R} |y|^{2n} \tilde{P}_t (dy)\leq c_n t^{2n/\alpha},\quad t>0.
\end{equation}
\end{lemat}

\dowod
Let $f_\varepsilon(\xi)=\int_{\R}(1-\cos(\scalp{y}{\xi}))\tilde{\nu}_\varepsilon(dy)$. For every $l\in\N$ and $k\in\{1,\dots,d\}$ we have
\begin{equation}\label{eq:Difff}
\frac{\partial^{2l}f_\varepsilon}{\partial\xi_k^{2l}}(0) = (-1)^{l+1} \int_{\R} y_k^{2l} \tilde{\nu}_\varepsilon(dy),
\end{equation}
and 
\begin{equation}\label{eq:Difff2}
\frac{\partial^{2l-1}f_\varepsilon}{\partial\xi_k^{2l-1}}(0) = 0.
\end{equation}
We use the Faa di Bruno's formula (see \cite{ConSav96}), (\ref{eq:Difff}) and (\ref{eq:Difff2}) 
to obtain
\begin{eqnarray}\label{eq:FaaDiB}
  \frac{\partial^{2n}}{\partial \xi_k^{2n}} \Fourier(\tilde{P}^{\varepsilon}_t)(0) & = &
   (2n)! \exp(-tf_\varepsilon(0))
    \sum_{j=1}^{2n} \sum_{\pi(2n,j)} \prod_{l=1}^{2n}\frac{\left(\frac{\partial^l (-tf_\varepsilon)}{\partial \xi_k^l}(0)\right)^{\lambda_l}}{(\lambda_l!)(l!)^{\lambda_l}}\nonumber \\
    &  =  & \sum_{j=1}^{2n} \sum_{\pi'(2n,j)}t^j\prod_{l=1}^{n}
    \frac{(-1)^n(2n)!}{(\lambda_{2l}!)((2l)!)^{\lambda_{2l}}} \left[\int_{\R}y_k^{2l}\tilde{\nu}_\varepsilon(dy)\right]^{\lambda_{2l}}
\end{eqnarray}
where 
$$
  \pi(2n,j)=\{(\lambda_1,\dots,\lambda_{2n}):\, \lambda_l\in \N_0,\,\sum_{l=1}^{2n}\lambda_l=j,\,\sum_{l=1}^{2n}l\lambda_l=2n\}
$$ and
$$
\pi'(2n,j)=\pi(2n,j)\cap\{(\lambda_1,\dots,\lambda_{2n}):\, \lambda_l\in \N_0,\,\lambda_1=\lambda_3=\dots=\lambda_{2n-1}=0\}.
$$
For every $n\in\N$ and $k\in\{1,\dots,d\}$ we have
\begin{eqnarray*}
  \int_{\R} y_k^{2n} \tilde{P}^{\varepsilon}_t (dy) & = &
  \frac{1}{i^{2n}}\left[\frac{\partial^{2n}}{\partial \xi_k^{2n}}  
  \Fourier(\tilde{P}^{\varepsilon}_t)(\xi)\right]_{\xi=0} \\
    &  =  &
  (-1)^n \left[\frac{\partial^{2n}}{\partial \xi_k^{2n}} \Fourier(\tilde{P}^{\varepsilon}_t)(\xi)\right]_{\xi=0},
\end{eqnarray*}
and this and (\ref{eq:FaaDiB}) yield
\begin{eqnarray}\label{eq:FaaDiB2}
	\int_{\R} |y|^{2n} \tilde{P}^{\varepsilon}_t (dy)
	& \leq & d^{n-1} \int_{\R} (y_1^{2n}+\dots+y_d^{2n})\tilde{P}^{\varepsilon}_t (dy) \nonumber \\ 
	&   =  & d^{n-1} \sum_{k=1}^d \sum_{j=1}^{2n} \sum_{\pi'(2n,j)}t^j\prod_{l=1}^{n}
    \frac{(2n)!}{(\lambda_{2l}!)((2l)!)^{\lambda_{2l}}} \left[\int_{\R}y_k^{2l}\tilde{\nu}_\varepsilon(dy)\right]^{\lambda_{2l}}.
	\end{eqnarray}
For every $k\in\{1,\dots,d\}$ and $l\in\N$ by (\ref{eq:AssAux1}) we have 
\begin{eqnarray*}
  \int_{\R} y_k^{2l}\tilde{\nu}_\varepsilon(dy)
  &  =   & \sum_{j=0}^{\infty} \int_{2^{-j-1}\varepsilon\leq |y| < 2^{-j}\varepsilon} 
           y_k^{2l} \tilde{\nu}_\varepsilon(dy) \\
  & \leq & \sum_{j=0}^{\infty} (2^{-j}\varepsilon)^{2l}\nu(B(0,2^{-j-1}\varepsilon)^c)\\
  & \leq & \varepsilon^{2l}\sum_{j=0}^{\infty} (2^{2l})^{-j}c(2^{-j-1}\varepsilon)^{-\alpha} = c \varepsilon^{2l-\alpha}.
\end{eqnarray*}
and by (\ref{eq:FaaDiB2}) we obtain
\begin{eqnarray*}
  \int_{\R} |y|^{2n} \tilde{P}^{\varepsilon}_t (dy)
  & \leq & \varepsilon^{2n} d^{n-1} \sum_{k=1}^d \sum_{j=1}^{2n}
           c(j,n)t^j\varepsilon^{-j\alpha}.
\end{eqnarray*}
Taking $\varepsilon=t^{1/\alpha}$ we get (\ref{eq:MomentsTildeP}).
\qed

\begin{lemat}\label{lm:TildepEstimate}
  If (\ref{eq:AssAux1}) is satisfied then for every $m\geq 1$ there exists a constant 
  $c_m=c(m)$ such that
  \begin{equation}\label{eq:TildepEstimate}
    \tilde{p}_t(y) \leq c(m)t^{-d/\alpha}(1+t^{-1/\alpha}|y|)^{-m}
    ,\quad y\in\R,\, t\in(0,1).
  \end{equation}
\end{lemat}
\dowod From (\ref{eq:AssAux1}) and (\ref{eq:FTildePEstimate}) we get
$\Fourier(\tilde{p}_t)(\xi)\leq c_1\exp(-c_2t|\xi|^\alpha)$ for $|\xi|>1$ and $\Fourier(\tilde{p}_t)(\xi)\leq 1$
for $|\xi|\leq 1$.
We denote $g_t(y)=t^{d/\alpha}\tilde{p}_t(t^{1/\alpha}y)$. For
every $j\in\{1,\dots,d\}$ we obtain
\begin{eqnarray}\label{eq:gprime}
  \left|\frac{\partial g_t}{\partial y_j}(y)\right|
  &   =  & \left|t^{d/\alpha}(2\pi)^{-d}\int_{\R}
           (-i)t^{1/\alpha}\xi_j e^{-it^{1/\alpha}\scalp{y}{\xi}}
           \Fourier(\tilde{p}_t)(\xi)d\xi\right| \nonumber \\
  & \leq & c t^{\frac{1+d}{\alpha}}\left( \int_{|\xi|\leq 1}|\xi_j|d\xi
             +\int_{|\xi|>1}|\xi_j|e^{-c_2 t|\xi|^\alpha}d\xi\right) \nonumber \\
  &   =  & c t^{\frac{1+d}{\alpha}}
           +c\int_{|u|>t^{1/\alpha}}|u_j|e^{-c_2 |u|^\alpha}du \nonumber  \\
  & \leq & c,\quad y\in\R,\, t\in(0,1).
\end{eqnarray}
Similarly we get $g_t(y)\leq c$ for $y\in\R$, $t\in(0,1)$.
By Lemma \ref{lm:MomentsTildeP} we obtain $\int_{\R} |y|^{2n} g_t(y)dy\leq c(n)$, and also
$\int_{\R} |y|^{2n-1} g_t(y)dy\leq \int_{|y|\leq 1} g_t(y)dy + \int_{|y|>1} |y|^{2n} g_t(y)dy\leq 1 +c(n)$
and by \cite[Lemma 9]{Pi} for every $m\geq 1$ we get
$$
  g_t(y) \leq c(m) (1+|y|)^{-m}, \quad y\in\R,\, t\in(0,1),
$$
which clearly yields (\ref{eq:TildepEstimate}).
\qed

We prove now the estimate of $\tilde{p}_t^\varepsilon(y)$ from below for small $y$, $\varepsilon=at^{1/\alpha}$
and $a\in(0,1]$.

\begin{lemat}\label{lm:LocalBelow}
If (\ref{eq:Ass7}) %and (\ref{eq:Ass8}) 
holds then there exist constants $c$, $c_1$, $c_2$ such that
\begin{equation}\label{eq:TildaBelow}
  \tilde{p}^{at^{1/\alpha}}_t(y) \geq c t^{-d/\alpha},  
\end{equation}
provided $|y|\leq c_1e^{-c_2 a^{-\alpha}}t^{1/\alpha}$, $t\in (0,1)$, and $a\in (0,1]$.  
\end{lemat}
\dowod It follows from (\ref{eq:Ass7}) that
\begin{eqnarray}\label{eq:Ass8}
  \Phi(\xi)  & \leq & \frac{1}{2}\int_{|y|\leq 1/|\xi|}|\scalp{\xi}{y}|^2 \nu(dy) +2\nu(B(0,1/|\xi|)^c) \nonumber \\
             & \leq & \frac{1}{2}|\xi|^2 \sum_{j=0}^{\infty} \int_{\frac{2^{-j-1}}{|\xi|}<|y|\leq \frac{2^{-j}}{|\xi|}} |y|^2 \nu(dy) + 2\nu(B(0,1/|\xi|)^c) \nonumber \\
             & \leq & \frac{1}{2}|\xi|^2 \sum_{j=0}^{\infty} 2^{-2j}|\xi|^{-2}\nu(B(0,\frac{2^{-j-1}}{|\xi|})^c)+2\nu(B(0,1/|\xi|)^c)
                     \nonumber \\
             & \leq &  c |\xi|^{\alpha}, \quad |\xi|>1.
\end{eqnarray}
Let $g_t(y)=t^{d/\alpha}\tilde{p}^{at^{1/\alpha}}_t (t^{1/\alpha}y)$. We have
$$
  \Fourier(\tilde{p}^{at^{1/\alpha}}_t)(\xi)
   \geq  \Fourier(p_t)(\xi),\quad \xi\in\R,\, t>0,
$$ and this and 
(\ref{eq:Ass8}) 
yield
\begin{eqnarray}\label{eq:g0}
  g_t(0) = t^{d/\alpha}\tilde{p}^{at^{1/\alpha}}_t(0) 
         & \geq & t^{d/\alpha} p_t(0) \nonumber \\
         & = & t^{d/\alpha}(2\pi)^{-d} \int_{\R} e^{-t\Phi(\xi)}d\xi \nonumber \\
         & \geq & t^{d/\alpha}(2\pi)^{-d}  \int_{|\xi|\geq 1}e^{-ct|\xi|^\alpha}d\xi \nonumber \\
         & \geq & c_0 >0, \quad t\in(0,1).
\end{eqnarray}
For every $j\in\{1,\dots,d\}$ by (\ref{eq:FTildePEstimate}) and (\ref{eq:Ass7}) we get	
\begin{eqnarray*}
  \left|\frac{\partial g_t}{\partial y_j}(y)\right|
  &   =  & \left|t^{d/\alpha}(2\pi)^{-d}\int_{\R}
           (-i)t^{1/\alpha}\xi_j e^{-it^{1/\alpha}\scalp{y}{\xi}}
           \Fourier(\tilde{p}^{at^{1/\alpha}}_t)(\xi)d\xi\right| \\
  & \leq & c_3 e^{c_2a^{-\alpha}}t^{\frac{1+d}{\alpha}}\left( \int_{|\xi|\leq 1}|\xi_j|d\xi
             +\int_{|\xi|>1}|\xi_j|e^{-c_4 t|\xi|^\alpha}d\xi\right) \\
  & \leq & c_5 e^{c_2a^{-\alpha}},\quad y\in\R,\, t\in(0,1).
\end{eqnarray*}
It follows that
$$
  g_t(y)\geq c_0-dc_5e^{c_2a^{-\alpha}}|y|\geq c_0/2,
$$
provided  $|y|\leq c_1e^{-c_2a^{-\alpha}}$, for some $c_1$
which clearly yields (\ref{eq:TildaBelow}).
\qed

%and by (\ref{eq:g0}) and (\ref{eq:gprime}) we get (\ref{eq:TildaBelow}).

%--------------------------------------------------------------------
\section{Bounded part}\label{sc:BP}

For a measure $\lambda$ on $\R$, $|\lambda|$ denotes
its total mass.
When $|\lambda|<\infty$ and $n=1,2,\ldots$ we let $\lambda^{n*}$
denote the $n$-fold convolution of $\lambda$ with itself:
$$
  \lambda^{n*}(f)=\int f(x_1+x_2+\dots+x_n)\lambda(dx_1)\lambda(dx_2)\ldots\lambda(dx_n)\,.
$$
We also let $\lambda^{0*}=\delta_0$, the evaluation at $0$.

For $\varepsilon>0$ we denote $\nubounded{\varepsilon}=\indyk{B(0,\varepsilon)^c}\nu$. 
We consider the corresponding semigroup of measures $\{\bar{P}^\varepsilon_t,\; t\geq 0\}$ such that
\begin{equation}\label{eq:FPbar}
  \Fourier(\bar{P}^\varepsilon_t)(\xi) =
  \exp\left(t \int (\cos(\scalp{\xi}{y})-1)
  \nubounded{\varepsilon}(dy)\right)\, ,
  \quad \xi\in\R\, .
\end{equation}
Note that
\begin{eqnarray}\label{eq:exp}
  \bar{P}^\varepsilon_t 
  &=& \exp(t(\bar{\nu_\varepsilon}-|\nubounded{\varepsilon}|\delta_0))
      =  \sum_{n=0}^\infty
\frac{t^n\left(\bar{\nu_\varepsilon}-|\nubounded{\varepsilon}|\delta_0)\right)^{n*}}{n!}\\
&=& e^{-t|\nubounded{\varepsilon}|}
    \sum_{n=0}^\infty
    \frac{t^n\nubounded{\varepsilon}^{n*}}{n!}\,,\quad t\geq 0\, , \nonumber
\end{eqnarray}
and
\begin{displaymath}%\label{erozklad}
P_t=\tilde{P}^\varepsilon_t \ast \bar{P}^\varepsilon_t\,,\quad t\geq 0.
\end{displaymath}
Of course,
\begin{equation}\label{e:wpt}
p_t=
\tilde{p}^\varepsilon_t
*
\bar{P}^\varepsilon_t
\,, \quad t>0.
\end{equation}

\begin{lemat}\label{lm:BoundedNuEst}
 If (\ref{eq:Ass1}), (\ref{eq:Ass2}) and (\ref{eq:Ass3}) are satisfied then there exists a constant $c$ such that
 \begin{equation}\label{eq:BoundedNuEst}
   \nubounded{\varepsilon}^{n*}(B(x,r)) 
    \leq  c^n r^\gamma (\varepsilon^{-\alpha}\qu(\varepsilon))^{n-1}|x|^{-\alpha-\gamma}\qu(|x|),
 \end{equation}
provided $x\in\R\setminus\{0\}$, $\varepsilon>0$, $n\in\N$, and $r\leq\max(\frac{\varepsilon}{3},\frac{|x|}{5^n} )$.
\end{lemat}
\dowod We proceed by induction. Note that (\ref{eq:BoundedNuEst}) for $n=1$ holds by (\ref{eq:Col1})
and (\ref{eq:Col2}). 
Let $c_0$ and $n$ be such that (\ref{eq:BoundedNuEst}) is satisfied with $c=c_0$. We first assume that
$r\leq\varepsilon/3$. By (\ref{eq:Col1}) we have
\begin{eqnarray*}
  \nubounded{\varepsilon}^{(n+1)*}(B(x,r))
  &   =  & \int\limits_{|x-y|>2\varepsilon/3} \nubounded{\varepsilon}(B(x-y,r))
           \nubounded{\varepsilon}^{n*}(dy) \\
  & \leq & \int\limits_{|x-y|>2\varepsilon/3} \nu(B(x-y,r))\nubounded{\varepsilon}^{n*}(dy) \\
  & \leq & c_1 r^\gamma \int\limits_{|x-y|>2\varepsilon/3}
            |x-y|^{-\alpha-\gamma}\qu(|x-y|)
           \nubounded{\varepsilon}^{n*}(dy)
\end{eqnarray*}
(note that $r<|x-y|/2$ provided $|x-y|>2\varepsilon/3$).
Now let $\varepsilon/3 < r \leq |x|/5^{n+1}$. Then $2r+\varepsilon<|x|/5^n$ and 
by induction and (\ref{eq:Col2})
\begin{eqnarray*}
  \int\limits_{|x-y|<2r+\varepsilon}\nubounded{\varepsilon}(B(x-y,r))
  \nubounded{\varepsilon}^{n*}(dy)
  & \leq & |\nubounded{\varepsilon}| \nubounded{\varepsilon}^{n*}(B(x,2r+\varepsilon)) \\
  & \leq & c_2 \varepsilon^{-\alpha}\qu(\varepsilon)
           c_0^{n}(2r+\varepsilon)^\gamma (\varepsilon^{-\alpha}\qu(\varepsilon))^{n-1}|x|^{-\alpha-\gamma}\qu(|x|) \\
  & \leq & c_0^n c_3 r^{\gamma} (\varepsilon^{-\alpha}\qu(\varepsilon))^n|x|^{-\alpha-\gamma}\qu(|x|),
\end{eqnarray*}
for some $c_3$; and by (\ref{eq:Col1}) we get
\begin{eqnarray*}
  \int\limits_{|x-y|\geq 2r+\varepsilon}\nubounded{\varepsilon}(B(x-y,r))
  \nubounded{\varepsilon}^{n*}(dy)
  & \leq & \int\limits_{|x-y|\geq 2r+\varepsilon}\nu(B(x-y,r))
           \nubounded{\varepsilon}^{n*}(dy) \\
  & \leq & \int\limits_{|x-y|\geq 2r+\varepsilon} c_1 r^{\gamma}
           |x-y|^{-\alpha-\gamma}\qu(|x-y|) \nubounded{\varepsilon}^{n*}(dy) \\
  & \leq & c_1 r^{\gamma}\int\limits_{|x-y|>2\varepsilon/3} 
           |x-y|^{-\alpha-\gamma}\qu(|x-y|) \nubounded{\varepsilon}^{n*}(dy).
\end{eqnarray*}
From the above we have
\begin{eqnarray}\label{eq:Summary}
  \nubounded{\varepsilon}^{(n+1)*}(B(x,r)) 
  & \leq  & c_1 r^{\gamma}\int\limits_{|x-y|>2\varepsilon/3} 
            |x-y|^{-\alpha-\gamma} \qu(|x-y|)\nubounded{\varepsilon}^{n*}(dy) \\
  &       &  + c_0^n c_3 r^{\gamma} (\varepsilon^{-\alpha}\qu(\varepsilon))^n
           |x|^{-\alpha-\gamma}\qu(|x|),\nonumber
\end{eqnarray}
for all $0<r\leq\max(\varepsilon/3,|x|/5^{n+1})$.

Let $L_\varepsilon=\czc{\log_5(\frac{3|x|}{2\varepsilon})}$. If $2\varepsilon/3<|x|/5^n$ then we get by (\ref{eq:Ass2a}) and induction
$$
\begin{array}{lcl}
  &       \int\limits_{2\varepsilon/3<|x-y|<|x|/5^n}&
  |x-y|^{-\alpha-\gamma}\qu(|x-y|) \nubounded{\varepsilon}^{n*}(dy)    \\
  & \leq & \qu(2\varepsilon/3)\sum\limits_{k=n}^{L_\varepsilon}\;
           \int\limits_{|x|/5^{k+1}\leq|x-y|<|x|/5^k}
           |x-y|^{-\alpha-\gamma} \nubounded{\varepsilon}^{n*}(dy)   \\
  & \leq & \qu(2\varepsilon/3) \sum\limits_{k=n}^{L_\varepsilon}
           (5^{k+1})^{\alpha+\gamma}|x|^{-\alpha-\gamma} 
           \nubounded{\varepsilon}^{n*}(B(x,|x|/5^k)) \\
  & \leq & c_0^n \qu(2\varepsilon/3) 5^{\alpha+\gamma} 
           (\varepsilon^{-\alpha}\qu(\varepsilon))^{n-1}|x|^{-2\alpha-\gamma} \qu(|x|)
           \sum\limits_{k=1}^{L_\varepsilon} 5^{k\alpha} \\
  & \leq & c_0^n c_4 (\varepsilon^{-\alpha}\qu(\varepsilon))^n|x|^{-\alpha-\gamma}\qu(|x|)\,.
\end{array}
$$
Also, by (\ref{eq:Ass2a}) and (\ref{eq:Col2})
$$
\begin{array}{lcl}
  \int\limits_{|x-y|\geq |x|/5^n} |x-y|^{-\alpha-\gamma}\qu(|x-y|)\nubounded{\varepsilon}^{n*}(dy)
  & \leq & (5^{\alpha+\gamma})^n |x|^{-\alpha-\gamma}c_5\qu(|x|)5^{n\eta}|\nubounded{\varepsilon}^{n*}| \\
  & \leq & c_5 c_0^n (\varepsilon^{-\alpha}\qu(\varepsilon))^n |x|^{-\alpha-\gamma}\qu(|x|),
\end{array}
$$
by taking large $c_0$.
We get
$$
  \int\limits_{|x-y|>2\varepsilon/3}
  |x-y|^{-\alpha-\gamma} \qu(|x-y|)\nubounded{\varepsilon}^{n*}(dy) \leq c_0^n (\varepsilon^{-\alpha}\qu(\varepsilon))^n (c_4+c_5)|x|^{-\alpha-\gamma}\qu(|x|),
$$
and (\ref{eq:Summary}) yields
$$
  \nubounded{\varepsilon}^{(n+1)*}(B(x,r)) \leq c_0^{n+1} r^{\gamma}(\varepsilon^{-\alpha}\qu(\varepsilon))^n|x|^{-\alpha-\gamma}\qu(|x|).
$$ 
\qed

\begin{wniosek}\label{co:BoundedNuEst}
If (\ref{eq:Ass1}), (\ref{eq:Ass2}) and (\ref{eq:Ass3}) are satisfied then there exists $c$ such that
  \begin{eqnarray}\label{eq:BoundedNuSumm}
    \nubounded{\varepsilon}^{n*}(B(x,r)) 
    & \leq & c^n (\varepsilon^{-\alpha}\qu(\varepsilon))^{n-1} r^\gamma(1+\frac{\varepsilon^{-\alpha}\qu(\varepsilon)}{r^{-\alpha}\qu(r)})
             |x|^{-\alpha-\gamma}\qu(|x|), \\
    &      & \quad x\in\R,\,r>0,\,\varepsilon>0,\,n\in\N.\nonumber
  \end{eqnarray}
\end{wniosek}
\dowod
If $r\leq |x|/5^n$ then (\ref{eq:BoundedNuSumm}) follows directly from Lemma \ref{lm:BoundedNuEst}.
If $r > |x|/5^n$ then by (\ref{eq:Col2}) and (\ref{eq:Ass2a}) we have
\begin{eqnarray*}
  \nubounded{\varepsilon}^{n*}(B(x,r)) 
  & \leq &  |\nubounded{\varepsilon}^{n*}|\leq c^n (\varepsilon^{-\alpha}\qu(\varepsilon))^n \\
  & \leq & c^n (\varepsilon^{-\alpha}\qu(\varepsilon))^n
           \left(\frac{r5^n}{|x|}\right)^{\gamma+\alpha}\frac{\qu(|x|)5^{n\eta}}{\qu(r)}.
\end{eqnarray*}
\qed

We will denote $$\bar{P}_t=\bar{P}_t^{t^{1/\alpha}}.$$

\begin{wniosek}\label{co:BoundedPEst}
If (\ref{eq:Ass1}), (\ref{eq:Ass2}) and (\ref{eq:Ass3}) are satisfied then there exists $c$ such that
$$
  \bar{P}_t(B(x,r)) \leq
  c \frac{t}{\qu(t^{1/\alpha})} r^{\gamma}\left(1+\frac{t^{-1}\qu(t^{1/\alpha})}{r^{-\alpha}\qu(r)}\right)|x|^{-\alpha-\gamma}\qu(|x|)
  \,,\quad  r>0\,,\; t>0\,,\; x\in\R.
$$
\end{wniosek}
\dowod Corollary \ref{co:BoundedNuEst} and (\ref{eq:exp}) yield
\begin{eqnarray*}
  \bar{P}_t(B(x,r))
  & \leq & \sum\limits_{n=0}^\infty
           \frac{t^n \nubounded{t^{1/\alpha}}^{n*}(B(x,r))}{n!} \\
  & \leq & \sum\limits_{n=0}^\infty
           \frac{c^n t (\qu(t^{1/\alpha}))^{n-1} r^\gamma\left(1+\frac{t^{-1}\qu(t^{1/\alpha})}{r^{-\alpha}\qu(r)}\right)|x|^{-\alpha-\gamma}\qu(|x|)}{n!} \\
  & \leq & e^{c\|q\|_\infty} \frac{t}{\qu(t^{1/\alpha})} r^\gamma \left(1+\frac{t^{-1}\qu(t^{1/\alpha})}{r^{-\alpha}\qu(r)}\right)|x|^{-\alpha-\gamma}\qu(|x|).
\end{eqnarray*}
\qed

%----------------------------------------------------------------
\section{Proofs of main results}\label{sc:MainR}

In the following two lemmas we investigate the transition densities at the origin.
The behaviour of $\Phi$ at infinity determinates the asymptotic of $p_t$ in small time
whereas considering the rate of decay of $\Phi$ at the origin we obtain the estimates of $p_t$ for large $t$.

\begin{lemat}\label{lm:pBasicEstimate} If (\ref{eq:Ass0}) holds then there exists $c$ such that
  $$
    p_t(x) \leq c t^{-d/\alpha},\quad x\in\R,\, t\in (0,1).
  $$
If (\ref{eq:Ass7}) 
%and (\ref{eq:Ass8}) 
is satisfied then there exist constants $c$ and $\delta$ such that
\begin{equation}\label{eq:ptb}
  p_t(x) \geq c t^{-d/\alpha}, \quad |x| \leq \delta t^{1/\alpha},\, t\in (0,1).
\end{equation}
\end{lemat}
\dowod 
By (\ref{eq:Ass0}) we have
\begin{eqnarray*}
	p_t(x) &  =   & (2\pi)^{-d} \int_{\R} e^{-i\scalp{x}{\xi}}e^{-t\Phi(\xi)}d\xi \\
	       & \leq & (2\pi)^{-d} \int_{\R} e^{-t\Phi(\xi)}d\xi \\
	       & \leq & (2\pi)^{-d} \left( \int_{|\xi|\leq 1} d\xi 
	                                  + \int_{|\xi|\geq 1} e^{-c_1 t|\xi|^\alpha}d\xi \right) \\
	       &   =  & (2\pi)^{-d} \left( c_2
                    + t^{-d/\alpha} \int_{|u|\geq t^{1/\alpha}}e^{-c_1|u|^\alpha}du\right)\\
         & \leq & c_3 \left( 1 + t^{-d/\alpha}\right)\\
         & \leq & 2c_3 t^{-d/\alpha} ,\quad x\in\R,\, t\in (0,1).
\end{eqnarray*}
%Let $g_t(x)=t^{d/\alpha}p_t(t^{1/\alpha}x)$. For every $t\in(0,1)$ the function $g_t$ is a density
%of the infinitely divisible distribution with the L\'evy measure $\mu_t(A)=t\nu(t^{1/\alpha}A)$.
%By (\ref{eq:Ass8}) we have (cf. (\ref{eq:g0}))
%$$
% g_t(0) \geq \hat{c} >0,\quad t\in (0,1).
%$$
%For every $j\in\{1,\dots,d\}$ by (\ref{eq:Ass0}) we obtain
%\begin{eqnarray*}
%  \left|\frac{\partial g_t}{\partial y_j}(y)\right|
%  &   =  & \left|t^{d/\alpha}(2\pi)^{-d}\int_{\R}
%           (-i)\xi_jt^{1/\alpha}e^{-it^{1/\alpha}\scalp{y}{\xi}}
%           \Fourier(p_t)(\xi)d\xi\right| \\
%  & \leq & c t^{\frac{1+d}{\alpha}}\left( \int_{|\xi|<1}|\xi_j|d\xi
%             +\int_{|\xi|>1}|\xi_j|e^{-ct|\xi|^\alpha}d\xi\right) \\
%  &   =  & c t^{\frac{1+d}{\alpha}}
%           +c\int_{\R}|\xi_j|e^{-c|\xi|^\alpha}d\xi \\
%  & \leq & c,\quad t\in(0,1),
%\end{eqnarray*}
%and we get
%$$
%  g_t(y)\geq \hat{c}/2,\quad |y|<\delta,\, t\in(0,1),
%$$
%for some constant $\delta>0$. This yields (\ref{eq:ptb}).
The proof of (\ref{eq:ptb}) is analogous to the proof of Lemma \ref{lm:LocalBelow} and
therefore we omit
the details.
\qed

\begin{lemat}\label{lm:pBasicEstimateBigt} If (\ref{eq:Ass4}) holds then there exists $c$ such that
  $$
    p_t(x) \leq c t^{-d/\beta},\quad x\in\R,\, t\in (1,\infty).
  $$
If furthermore (\ref{eq:Ass10}) is satisfied then there exist constants $c$ and $\delta$ such that
\begin{equation}\label{eq:pBigtBasisbelow}
  p_t(x) \geq c t^{-d/\beta}, \quad |x| \leq \delta t^{1/\beta},\, t\in (1,\infty).
\end{equation}
\end{lemat}
\dowod By (\ref{eq:Ass0}) and (\ref{eq:Ass4}) we have
\begin{eqnarray*}
	p_t(x) &  =   & (2\pi)^{-d} \int_{\R} e^{-i\scalp{x}{\xi}}e^{-t\Phi(\xi)}d\xi \\
	       & \leq & (2\pi)^{-d} \int_{\R} e^{-t\Phi(\xi)}d\xi \\
	       & \leq & (2\pi)^{-d} \left( \int_{|\xi|\leq 1} e^{-ct|\xi|^\beta} d\xi 
	                                  + \int_{|\xi|\geq 1} e^{-ct|\xi|^\alpha}d\xi \right) \\
	       &   =  & (2\pi)^{-d} \left( t^{-d/\beta} \int_{|u|\leq t^{1/\beta}}e^{-c|u|^\beta}du
                    + t^{-d/\alpha} \int_{|u|\geq t^{1/\alpha}}e^{-c|u|^\alpha}du\right)\\
         & \leq & c \left( t^{-d/\beta} + t^{-d/\alpha}\right)\\
         & \leq & c t^{-d/\beta} ,\quad x\in\R,\, t\in (1,\infty).
\end{eqnarray*}
We omit the proof of (\ref{eq:pBigtBasisbelow}) since it is analogous to the proof of Lemma \ref{lm:LocalBelow}.
\qed

%This yields that there exists a constant $\delta>0$ such that 
%$ g_t(x)\geq c/2 >0$ for $|x|\leq \delta$ and the Lemma follows.
%\qed

{\it Proof of the Theorem \ref{th:PEstimates}.}
Similarly like in the proof of \cite[Theorem 3]{Pi} and
\cite[Lemma 6]{BS2007} for $m=\gamma+\alpha+\eta+1$ by Lemma \ref{lm:TildepEstimate}, 
Corollary \ref{co:BoundedPEst} and (\ref{eq:Ass2a}) we obtain
\begin{eqnarray*}
  p_t(x)  &  =   &  \int_{\R} \tilde{p}_t(x-z)\bar{P}_t(dz) \\
                   & \leq & c \int_{\R} t^{-d/\alpha}(1+t^{-1/\alpha}|x-z|)^{-m} 
                              \bar{P}_t(dz)\\
                   &  =   & c t^{-d/\alpha} \int_0^1 
                            \bar{P}_t(\{z\,:\; (1+t^{-1/\alpha}|x-z|)^{-m}>s\})ds \\
                   & \leq & c t^{-d/\alpha} \int_0^1 \bar{P}_t(B(x, t^{1/\alpha}s^{-1/m}))ds\\
                   & \leq & c t^{-d/\alpha} \int_0^{1}\frac{t^{1+\gamma/\alpha}}{\qu(t^{1/\alpha})}s^{-\gamma/m}\left(1+\frac{s^{-\alpha/m}\qu(t^{1/\alpha})}{\qu(t^{1/\alpha}s^{-1/m})}\right)
                            |x|^{-\alpha-\gamma}\qu(|x|)ds \\
                   & \leq & c \frac{t^{1+\frac{\gamma-d}{\alpha}}}{\qu(1)}|x|^{-\alpha-\gamma}\qu(|x|)
                            \left[\int_0^1 s^{-\gamma/m}ds +
                                  \int_0^1 s^{-(\gamma+\alpha+\eta)/m}ds
                            \right] \\
                   &  =   & ct^{1+\frac{\gamma-d}{\alpha}}|x|^{-\alpha-\gamma}\qu(|x|), \quad x\in\R,\, t\in(0,1).
\end{eqnarray*}
This, together with Lemma \ref{lm:pBasicEstimate} gives (\ref{eq:AbDensityEstimate}).
\qed \\

\begin{lemat}\label{lm:Propnu}
  If (\ref{eq:Ass1}), (\ref{eq:Ass2}), (\ref{eq:Ass3}) and (\ref{eq:Ass5}) hold then there exists a constant c such that
  \begin{equation}\label{eq:Col3}
    \qu(s) \leq c s^{\alpha-\beta},\quad s\geq 1,
  \end{equation}
and
  \begin{equation}\label{eq:Col4}
    \int_{|y|<r} |y|^2 \nu(dy) \leq c r^{2-\beta},\quad r\geq 1.
  \end{equation}
\end{lemat}
\dowod If $\beta=\alpha$ then (\ref{eq:Col3}) follows from the boundedness of $\qu$. For $\beta>\alpha$ 
by (\ref{eq:Ass5}) we have
\begin{eqnarray*}
  \infty > \int_1^\infty u^{\beta-\alpha-1}\qu(u) du
  & \geq & \int_1^s u^{\beta-\alpha-1}\qu(u) du \\
  & \geq & \qu(s) \int_1^s u^{\beta-\alpha-1} du \\
  & \geq & c \qu(s)s^{\beta-\alpha},\quad s\geq 2.
\end{eqnarray*}
This and the boundedness of $\qu$ yield (\ref{eq:Col3}).\\
By (\ref{eq:Col2}) we have
\begin{eqnarray*}
  \int_{|y|<r}|y|^2 \nu(dy)
  & \leq & \int_{|y|<1}|y|^2 \nu(dy) + 
           \sum_{j=0}^{\czc{\log_2 r}} \int_{2^j \leq |y|<2^{j+1}}|y|^2 \nu(dy) \\
  & \leq & c_1 +  \sum_{j=0}^{\czc{\log_2 r}} 2^{2(j+1)}\nu(B(0,2^j)^c) \\
  & \leq & c_1 + c_2 \sum_{j=0}^{\czc{\log_2 r}} 2^{2j}2^{-j\alpha}\qu(2^{j}). \\
\end{eqnarray*}
If $\beta<2$ then (\ref{eq:Col4}) follows from (\ref{eq:Col3}). If $\beta=2$ then using
(\ref{eq:Ass2}) and (\ref{eq:Ass5}) we obtain
\begin{eqnarray*}
  \sum_{j=0}^{\czc{\log_2 r}} \left(2^{j}\right)^{2-\alpha}\qu(2^{j})
  & \leq & \sum_{j=0}^{\infty} \left(2^{j}\right)^{2-\alpha}\qu(2^{j}) \\
  & \leq & c \sum_{j=0}^{\infty} \int_{2^j}^{2^{j+1}} s^{1-\alpha}\qu(s)ds \\
  &   =  & c \int_1^\infty s^{1-\alpha}\qu(s) ds = c < \infty \\
\end{eqnarray*}
\qed

The proof of Theorem \ref{th:PEstimatesBigt} is similar to the proof of Theorem \ref{th:PEstimates} and
therefore we abbreviate below the details.

{\it Proof of the Theorem \ref{th:PEstimatesBigt}.} Let $\tilde{P}^*_t=\tilde{P}^{t^{1/\beta}}_t$, $\tilde{p}^*_t=\tilde{p}^{t^{1/\beta}}_t$, and $\bar{P}^*_t=\bar{P}_t^{t^{1/\beta}}$. By (\ref{eq:Col4}) for $l\in\N$ 
and $r>1$ we have
\begin{eqnarray*}
  \int_{\R} y_k^{2l}\tilde{\nu}_r(dy)
  & \leq &  \int_{|y|<r} |y|^{2l}\nu(dy) \\
  & \leq &  r^{2l-2} \int_{|y|<r} |y|^2 \nu(dy) \leq c r^{2l-\beta},
\end{eqnarray*}
and (\ref{eq:FaaDiB2}) yields
\begin{eqnarray*}
  \int_{\R} |y|^{2n} \tilde{P}^{r}_t (dy)
  & \leq & r^{2n} d^{n-1} \sum_{k=1}^d \sum_{j=1}^{2n}
           c(j,n)t^j r^{-j\beta}.
\end{eqnarray*}
Taking $r=t^{1/\beta}$ we get
\begin{equation}\label{eq:PStarMoment}
  \int_{\R} |y|^{2n} \tilde{P}^*_t (dy)
   \leq  c(n)t^{2n/\beta},\quad t\geq 1.
\end{equation}
From (\ref{eq:Col2}) and (\ref{eq:Col3}) we obtain $\nu(B(0,r)^c)\leq c r^{-\beta}$ for $r\geq 1$, and
similarly like in the proof of Lemma \ref{lm:TildepEstimate}, 
using (\ref{eq:Ass4}), (\ref{eq:FTildePEstimate}) and (\ref{eq:PStarMoment}), 
for every $m\geq 1$ we get
$$
   \tilde{p}^*_t(y) \leq c(m)t^{-d/\beta}(1+t^{-1/\beta}|y|)^{-m}
    ,\quad y\in\R,\, t\in(1,\infty).
$$
By Corollary \ref{co:BoundedNuEst} and (\ref{eq:Col3}) we get
$$
  \bar{P}^*_t(B(x,r)) \leq
  c t r^{\gamma}\left(1+\frac{t^{-\alpha/\beta}\qu(t^{1/\beta})}{r^{-\alpha}\qu(r)}\right)|x|^{-\alpha-\gamma}\qu(|x|)
  \,,\quad  x\in\R,\, t\geq 1,\, r>0,
$$
and for $m=\gamma+\alpha+\eta+1$ we obtain
\begin{eqnarray*}
  p_t(x)  &  =   &  \int_{\R} \tilde{p}^*_t(x-z)\bar{P}^*_t(dz) \\
                   & \leq & c \int_{\R} t^{-d/\beta}(1+t^{-1/\beta}|x-z|)^{-m} 
                              \bar{P}^*_t(dz)\\
                   %&  =   & c t^{-d/\beta} \int_0^1 
                   %         \bar{P}^*_t(\{z\,:\; (1+t^{-1/\beta}|z-y|)^{-m}>s\})ds \\
                   & \leq & c t^{-d/\beta} \int_0^1 \bar{P}^*_t(B(x, t^{1/\beta}s^{-1/m}))ds\\
                   & \leq & c t^{-d/\beta} \int_0^{1}t^{1+\gamma/\beta}s^{-\gamma/m}\left(1+\frac{s^{-\alpha/m}\qu(t^{1/\beta})}{\qu(t^{1/\beta}s^{-1/m})}\right)
                            |x|^{-\alpha-\gamma}\qu(|x|)ds \\
                   &   \leq   & c t^{1+\frac{\gamma-d}{\beta}} |x|^{-\alpha-\gamma}\qu(|x|)
                            \left[\int_0^1 s^{-\gamma/m}ds +
                                  c\int_0^1 s^{-(\gamma+\alpha+\eta)/m}ds
                            \right] \\
                   &  =   & ct^{1+\frac{\gamma-d}{\beta}}|x|^{-\alpha-\gamma}\qu(|x|), \quad x\in\R,\, t\in(1,\infty).
\end{eqnarray*}
We get (\ref{eq:AbDensityEstimateBigt}) by this and the Lemma \ref{lm:pBasicEstimateBigt}.
\qed

{\it Proof of the Theorem \ref{th:PEstimatesB}.}
Let $a\in(0,1)$ and $t\in(0,1)$.
For $|x|>r+at^{1/\alpha}$ by (\ref{eq:exp}) and (\ref{eq:Ass7}) we get
\begin{equation}\label{eq:BarNuBelow}
  \bar{P}^{at^{1/\alpha}}_t(B(x,r)) \geq e^{-ca^{-\alpha}} t \bar{\nu}_{at^{1/\alpha}}(B(x,r)) = e^{-ca^{-\alpha}} t \nu (B(x,r)).
\end{equation}
This, Lemma \ref{lm:LocalBelow} and (\ref{eq:Ass6}) for $x\in A$ yield 
\begin{eqnarray*}
  p_t(x) 
  &  =    & \tilde{p}^{at^{1/\alpha}}_t * \bar{P}_t^{at^{1/\alpha}}(x) \\
  &  =   & \int \tilde{p}_t^{at^{1/\alpha}}(x-z)\bar{P}_t^{at^{1/\alpha}}(dz) \\
  & \geq & c \int_{|z-x|<c_1e^{-c_2a^{-\alpha}}t^{1/\alpha}} t^{-d/\alpha} \bar{P}_t^{at^{1/\alpha}}(dz) \\
  &  =   & c t^{-d/\alpha}\bar{P}_t^{at^{1/\alpha}}(B(x,c_1e^{-c_2a^{-\alpha}}t^{1/\alpha})) \\
  & \geq & c(a) t^{1+\frac{\gamma-d}{\alpha}}|x|^{-\alpha-\gamma}\qd(|x|), 
\end{eqnarray*}
provided $|x|>(a+c_1e^{-c_2a^{-\alpha}})t^{1/\alpha}$.
By Lemma \ref{lm:pBasicEstimate} we have 
$p_t(x)\geq c t^{-d/\alpha}$ for $|x|<\delta t^{1/\alpha}$. We choose $a\in(0,1)$ such that
$a+c_1e^{-c_2a^{-\alpha}}\leq\delta$ and we obtain (\ref{eq:BeDensityEstimate}).
\qed

{\it Proof the Theorem \ref{th:PEstimatesBBigt}.}
If follows from (\ref{eq:Ass9}) and (\ref{eq:Ass10}) that 
$\tilde{p}^{at^{1/\beta}}_t (y)\geq c t^{-d/\beta}$ for $t \geq 1$, $a\in(0,1)$, and $|y|\leq c_1 e^{-c_2a^{-\beta}}t^{1/\beta}$ (cf. Lemma \ref{lm:LocalBelow}). Using this and Lemma \ref{lm:pBasicEstimateBigt} we get (\ref{eq:BeDensityEstimateBigt}) in the same way as (\ref{eq:BeDensityEstimate}) above.
\qed

%------------------------------------------------------------------------------------------------------
\section{Tempered stable processes}\label{sec:Temp}

In what follows we assume that
\begin{equation}\label{eq:Ex1}
  \nu(D) = \int_{\sfera} \int_0^\infty  \indyk{D}(s\theta) s^{-1-\alpha}Q(\theta,s) ds \mu(d\theta),
  \quad D\subset\R,
\end{equation}
where $\mu$ is a bounded symmetric measure on $\sfera$ and $Q$ is a nonnegative bounded function such that $Q(-\theta,s)=Q(\theta,s)$, $\theta\in\sfera$, $s>0$. We consider the semigroup of probability measures $\{ P_t, t\geq 0 \}$ with the Fourier transform $\Fourier(P_t)(\xi)=\exp(-t\Phi(\xi))$, where
\begin{equation}
  \Phi(\xi) =  \int_{\R} \left(1-\cos(\scalp{\xi}{y})\right)\nu(dy),\quad \xi\in\R.
\end{equation}

In \cite{Ros07} J.~Rosi\'nski investigate the {\it tempered stable process} with the L\'evy measure given
by (\ref{eq:Ex1}) with $Q(\theta,\cdot)$ completely monotone for
every $\theta\in\sfera$. We do not assume here
that $Q$ is completely monotone and consider a wider class of processes.

\begin{lemat}\label{lm:Ex1PhiEst}
  If there exist $s_0\in (0,1)$ such that
\begin{equation}\label{eq:Ex1Ass1}
  \inf_{s\in(0,s_0),\eta\in\sfera} \int_{\sfera} |\scalp{\eta}{\theta}|^2 Q(\theta,s)\mu(d\theta) > 0,
\end{equation}
and
\begin{equation}\label{eq:Ex1Ass2}
  \int_{\R} |y|^2 \nu(dy) = \int_{\sfera} \int_0^\infty  s^{1-\alpha}Q(\theta,s)ds\mu(d\theta)<\infty,
\end{equation}
then
\begin{equation}\label{eq:Ex1PhiEst}
  \Phi(\xi) \approx |\xi|^{2}\wedge |\xi|^{\alpha},\quad \xi\in\R.
\end{equation}
\end{lemat}
\dowod
  By (\ref{eq:Ex1Ass2}) we have
$$
  \Phi(\xi)  \leq  \frac{1}{2}\int_{\R} |\scalp{\xi}{y}|^2 \nu(dy)  \leq  \frac{1}{2} |\xi|^2 \int_{\R} |y|^2 \nu(dy)     =   c |\xi|^2,\quad \xi\in\R.
$$
For $|\xi|>1$ by the boundedness of $Q$ and $\mu$ we obtain also
\begin{eqnarray*}
  \Phi(\xi)
  & \leq & \frac{1}{2} \int_{|y|\leq 1/|\xi|} |\scalp{\xi}{y}|^2 \nu(dy) + 2\int_{|y| > 1/|\xi|}\nu(dy)\\
  & \leq & \frac{1}{2} |\xi|^2 \int_0^{1/|\xi|}\int_{\sfera}s^{1-\alpha}Q(\theta,s)ds\mu(d\theta)+
           2 \int_{1/|\xi|}^\infty \int_{\sfera} s^{-1-\alpha} Q(\theta,s)ds\mu(d\theta) \\
  & \leq & c|\xi|^\alpha.
\end{eqnarray*}
Moreover by (\ref{eq:Ex1Ass1}) we get
\begin{eqnarray*}
  \Phi(\xi)
  & \geq & (1-\cos 1) \int_{|y|<1/|\xi|} |\scalp{\xi}{y}|^2 \nu(dy) \\
  &   =  & (1-\cos 1) \int_{\sfera} \int_0^{1/|\xi|}
           |s\scalp{\xi}{\theta}|^2 s^{-1-\alpha}Q(\theta,s)ds\mu(d\theta) \\
  & \geq & c |\xi|^2 \int_0^{\frac{1}{|\xi|}\wedge s_0}s^{1-\alpha}ds ,\quad \xi\in\R,
\end{eqnarray*}
and the Lemma follows.
\qed

It follows from Lemma~\ref{lm:Ex1PhiEst} that we can apply the Theorem~\ref{th:PEstimates} 
with
$\qu(s)\equiv const.$ to every L\'evy measure given by (\ref{eq:Ex1}) with $Q$
and $\mu$ satisfying (\ref{eq:Ass3}) and (\ref{eq:Ex1Ass1}).
In specific cases we can obtain certainly more precise results.

We call a measure $\lambda$ {\it degenerate} if there is a proper
linear subspace $M$ of $\R$ such that $\supp(\lambda)\subset M$; 
otherwise we call $\lambda$ {\it nondegenerate}. Note that if $\mu$
is nondegenerate and $Q(\theta,s)\geq c >0$ on $\supp(\mu)\times(0,s_0)$
for some $s_0>0$ then (\ref{eq:Ex1Ass1}) is satisfied.

The following theorem follows from Lemma~\ref{lm:Ex1PhiEst} and
the Theorems~\ref{th:PEstimates}, \ref{th:PEstimatesB}, \ref{th:PEstimatesBigt} and \ref{th:PEstimatesBBigt} with
$\beta=2$ and $\qu(s)=\qd(s)=(1+s)^{-m}$.
We omit the easy proof.

\begin{twierdzenie}
If there exists $m>2-\alpha$ such that
\begin{equation}\label{eq:Ex1Ass3}
  Q(\theta,s) \approx (1+s)^{-m},\quad s>0,\, \theta\in\sfera,
\end{equation}
and $c$ and $\gamma\in[1,d]$ such that
\begin{equation}\label{eq:Ex1Ass4}
  \mu(B(\theta,r)\cap\sfera) \leq c r^{\gamma-1},\quad r<1/2,\,\theta\in\sfera,
\end{equation}
and $\mu$ is nondegenerate  then
$$
  p_t(x) \leq C \min\left(t^{-d/\alpha},\frac{t^{1+\frac{\gamma-d}{\alpha}}}
                     {|x|^{\alpha+\gamma}(1+|x|)^m}\,\right), \quad x\in\R,\, t\in(0,1),
$$
and
$$
   p_t(x)\leq C \min\left(t^{-d/2},\frac{t^{1+\frac{\gamma-d}{2}}}
                     {|x|^{\alpha+\gamma}(1+|x|)^m}\,\right),\quad  
   x\in\R,\, t\in(1,\infty).
$$
If additionally there exist a constant $c$ and a set $A_0\subset\sfera$ such that
\begin{equation}\label{eq:Ex1Ass5}
  \mu(B(\theta,r)\cap\sfera) \geq c r^{\gamma-1},\quad \theta\in A_0,\, r<1/2,
\end{equation}
then
$$
  p_t(x) \geq C \min\left(t^{-d/\alpha},\frac{t^{1+\frac{\gamma-d}{\alpha}}}
                     {|x|^{\alpha+\gamma}(1+|x|)^m}\,\right), \quad x\in A,\, t\in(0,1),
$$
and
$$
   p_t(x)\geq C \min\left(t^{-d/2},\frac{t^{1+\frac{\gamma-d}{2}}}
                     {|x|^{\alpha+\gamma}(1+|x|)^m}\,\right),\quad  
   x\in A,\, t\in(1,\infty),
$$
where $A=\{r\theta:\: r>0,\, \theta\in A_0\}$.
\end{twierdzenie}

If $\mu$ is absolutely continuous with respect to the standard surface measure on $\sfera$ and
its density $g_{\mu}$ is such that $c^{-1}\leq g_{\mu}(\theta)\leq c$, $\theta\in\sfera$, for some constant
$c>0$ then we have $\mu(B(\theta,r))\approx r^{d-1}$, $\theta\in\sfera$, $r\leq 1/2$.
Therefore in this case for $Q(\theta,s)$ satisfying (\ref{eq:Ex1Ass3}) we obtain
$$
  p_t(x) \approx \min\left(t^{-d/\alpha},\frac{t}
                     {|x|^{\alpha+d}(1+|x|)^m}\,\right), \quad x\in\R,\, t\in(0,1),
$$
and
$$
   p_t(x) \approx \min\left(t^{-d/2},\frac{t}
                     {|x|^{\alpha+d}(1+|x|)^m}\,\right),\quad  
   x\in \R,\, t\in(1,\infty).
$$

We can also apply our results for $\nu$ given by (\ref{eq:Ex1}), nondegenerate $\mu$ satisfying (\ref{eq:Ex1Ass4}) and $Q$ such that
\begin{equation}\label{eq:Expo}
  c^{-1}(1+s)^a e^{-c_1 s} \leq Q(\theta,s) \leq c(1+s)^a e^{-c_2 s},\quad s>0,\,\theta\in\sfera,
\end{equation}
for some $a\geq 0$.
In this  case for every $m>0$ we also have
\begin{equation}
  Q(\theta,s) \leq c_m (1+s)^{-m},\quad s>0,\,\theta \in\sfera,
\end{equation}
and therefore from Lemma~\ref{lm:Ex1PhiEst} and
Theorems~\ref{th:PEstimates} and \ref{th:PEstimatesBigt} we obtain
\begin{equation}\label{eq:ExE1}
  p_t(x) \leq C_m \min\left(t^{-d/\alpha},\frac{t^{1+\frac{\gamma-d}{\alpha}}}
                     {|x|^{\alpha+\gamma}(1+|x|)^m}\,\right), \quad x\in\R,\, t\in(0,1),
\end{equation}
and
\begin{equation}\label{eq:ExE2}
   p_t(y)\leq  C_m \min\left(t^{-d/2},\frac{t^{1+\frac{\gamma-d}{2}}}
                     {|x|^{\alpha+\gamma}(1+|x|)^m}\,\right),\quad  
   x\in\R,\, t\in(1,\infty).
\end{equation}
for every $m>2-\alpha$. If we assume additionally (\ref{eq:Ex1Ass5}) then by Theorems~\ref{th:PEstimatesB} and \ref{th:PEstimatesBBigt} with $\qd(s)=(1+s)^a e^{-2c_2s}$ we get
\begin{equation}\label{eq:ExE3}
  p_t(x) \geq C \min\left(t^{-d/\alpha},\frac{t^{1+\frac{\gamma-d}{\alpha}}}
                     {|x|^{\alpha+\gamma}}\,(1+|x|)^ae^{-c_3 |x|}\right), \quad x\in A,\, t\in(0,1),
\end{equation}
and
\begin{equation}\label{eq:ExE4}
   p_t(x) \geq C \min\left(t^{-d/2},\frac{t^{1+\frac{\gamma-d}{2}}}
                     {|x|^{\alpha+\gamma}}\,(1+|x|)^ae^{-c_3 |x|}\right), \quad x\in A,\, t\in(1,\infty).
\end{equation}

We like to discuss the particular case of the relativistic $\alpha$--stable L\'evy process with
$\Phi_m(\xi)=(|\xi|^2+m^{2/\alpha})^{\alpha/2}-m$ for $m>0$ which is investigated, e.g., in \cite{Ryz}, \cite{KulczSiu} and \cite{GrzRyz08}.
We consider here only $m=1$ because
$$
  p_t^m(x)=m^{d/\alpha}p^1_{mt}(m^{1/\alpha}x),
$$
where $p^m_t$ denotes the transition densities corresponding to $\Phi_m$, see \cite{GrzRyz08}. The L\'evy measure in this case
has the form
\begin{eqnarray*}
  \nu(D) & =  & c_1 \int_{D} |y|^{-d-\alpha} K_{d,\alpha}(|y|)dy \\
         & =  & c_2 \int_{\sfera} \int_0^\infty  \indyk{D}(s\theta) s^{-1-\alpha} K_{d,\alpha}(s) ds \sigma(d\theta),
        \quad D\subset\R,
\end{eqnarray*}
where $\sigma$ is the standard isotropic surface measure on $\sfera$ and 
$$
  K_{d,\alpha}(s)=s^{d+\alpha}\int_0^\infty e^{-u}e^{-\frac{s^2}{4u}}u^{\frac{-2-d-\alpha}{2}}du,\quad s>0.
$$
We have (see \cite{GrzRyz08})
$$
 K_{d,\alpha}(s) \approx (1+s)^{\frac{d+\alpha-1}{2}}e^{-s}.
$$
This yields that (\ref{eq:Expo}) for $Q(\theta,s)=K_{d,\alpha}(s)$ and $a=\frac{d+\alpha-1}{2}$ holds 
and we obtain the estimates
(\ref{eq:ExE1}), (\ref{eq:ExE2}), (\ref{eq:ExE3}) and (\ref{eq:ExE4}) with $\gamma=d$, $A=\R$, and $a=\frac{d+\alpha-1}{2}$. The sharp estimates of the transition densities of the relativistic process are given
also in \cite{GrzMast}.

\textbf{Acknowledgement.}
The author is grateful to Prof. R. L. Schilling  for useful discussions and his hospitality during the stay in
Dresden and to Prof. K. Bogdan
for discussions and suggestions on the paper.

%\vspace{2mm}
%P. Sztonyk, Institute of Mathematics, Wroc\l{}aw University of Technology,
%Wybrze{\.z}e Wyspia\'nskiego 27, 50--370 Wroc\l{}aw, Poland

%{\it E-mail address: sztonyk@pwr.wroc.pl}

\end{document}